\documentclass[11pt]{article}
\usepackage{color}
\usepackage{graphicx}

\newfont{\frak}{eufm10 scaled\magstep1}
\newfont{\sfrak}{eufm8 scaled\magstep1}
\newfont{\bbb}{msbm10 scaled\magstephalf}

\newtheorem{thm}{Theorem}[section]
\newtheorem{prop}[thm]{Proposition}

\newtheorem{defn}[thm]{Definition}
\newtheorem{remark}[thm]{Remark}

\newcommand{\qer}{\nolinebreak\hfill{$\bigtriangledown$}\par\vspace{0.5\parskip}}

\newcounter{sect}

\def\R{\mbox{\bbb{R}}}
\def\Z{\mbox{\bbb{Z}}}

\def\V{\mbox{\bbb{V}}}
\def\RP{\mbox{\bbb{RP}}}

\def\gl{\mbox{\frak gl}}
\def\g{\mbox{\frak g}}

\def\z{\mbox{\frak z}}

\def\D{\mbox{\frak D}}
\def\z{\mbox{\frak z}}

\title{\sc The Cauchy problem for Lie-minimal surfaces}

\author{\sc Emilio Musso}

\begin{document}

\maketitle

\maketitle \footnotetext[1]{This research was partially supported
by the MIUR project \textit{Propriet\`a Geometriche delle
Variet\`a Reali e Complesse}, by the group GNSAGA of the INdAM,
and by the European Contract Human Potential Programme, Research
Training Network HPRN-CT-2000-00101 (EDGE).
}


\begin{abstract}
In the present paper we study the Lie sphere geometry of Legendre
surfaces by the method of moving frame and we prove an existence
theorem for real-analytic Lie-minimal Legendre surfaces.
\end{abstract}

\section*{Introduction} In his analysis \cite{Bl2} of Lie sphere geometry of
surfaces W.Blaschke proposed to study the variational problem for
the functional \begin{equation}\mathcal{B}:M\subset \R^3\to \int_M
\frac{\partial_1k_1\partial_2k_2}{(k_1-k_2)^2}du^1\wedge
du^2,\label{I1}\end{equation} on immersed surfaces $M\subset \R^3$
with no umbilical points, where $k_1$ and $k_2$ are the principal
curvatures and where $(u^1,u^2)$ are curvature line coordinates.
He also showed that the functional is invariant under Lie sphere
transformations. Recently, E.Ferapontov \cite{Fe1} reconsidered
this classical variational problem and showed that the critical
points of (\ref{I1}) do admit a spectral deformation. This work
was taken up by F.Burstall and U.Hertrich-Jeromin \cite{BH}, who
introduced a Lie-invariant Gauss map $\mathcal{D}:M\to \D$ with
values in the \it{Dupin manifold }\rm $\D$, that is the symmetric
space consisting of all $3$-dimensional subspaces of signature
$(2,1)$ in $\R^{(4,2)}$. They showed that $M\subset \R^3$ is a
critical point of the functional (\ref{I1}) if and only if its
Lie-invariant Gauss map is harmonic. This explains the origin of
the spectral deformation discovered by Ferapontov and suggests the
existence of a dressing action on the space of the critical points
of the Blaschke functional (see \cite{BH}). In this paper, we
study the Cauchy problem for Lie-minimal surfaces using the
invariance by Lie sphere transformations from the outset. From
this point of view the relevant objects of study are the Legendre
lifts in the space of contact elements rather than the surfaces
itself.

In $\S 1$ we consider the Kepler manifold \footnote{We adopt the
terminology introduced by J.M Souriau in \cite{Sou} and by
Guillelmin and Sternberg in \cite{GS}} $\mathcal{K}$, that is the
isotropic Grassmannian of the null-planes through
$0\in\R^{(4,2)}$, acted upon transitively by the Lie sphere group
$\widetilde{G}=SO(4,2)/\pm I$ of contact transformations. We also
consider the \it{Lie quadric }\rm $\mathcal{Q}$ (i.e. the
projectivization of the light cone $\mathcal{L}$ of $\R^{(4,2)}$)
and the \it{Dupin manifold }\rm $\D$. We apply moving frame to
study Legendre surfaces in the Kepler manifold. Any Legendre
surface $M$ can be parameterized by a pair of mappings
$\phi_0,\phi_1:M\to \mathcal{Q}$ satisfying $\langle
\phi_0,\phi_1\rangle=0$ and $\langle d\phi_0,\phi_1\rangle = 0$.
Most of the Lie-invariant properties of $M$ are determined by the
sheaf $\mathcal{S}$ of quadratic forms spanned by $\langle
d\phi_0,d\phi_0\rangle$ and by $\langle d\phi_1,d\phi_1\rangle$.
If the stalks $\mathcal{S}_m$ are $2$-dimensional, for every $m\in
M$ then, the tautological bundle $$\mathcal{U}(M)=\{(\ell,V)\in
M\times \R^{(4,2)}:V\in \ell\}\to M$$ has a natural splitting into
the direct sum of two line sub-bundles $\Sigma_0(M)$ and
$\Sigma_1(M)$. The maps $\sigma_0,\sigma_1:M\to \mathcal{Q}$
induced by $\Sigma_0(M)$ and by $\Sigma_1(M)$ are the two
\it{curvature sphere mappings }\rm of the surface. If $\sigma_0$
and $\sigma_1$ are everywhere of maximal rank, then $M$ is said to
be \it{non-degenerate}\rm. Geometrically, this condition means
that the Euclidean projection of $M$ can not be obtained as the
envelope of a $1$-parameter family of spheres. In $\S 2$ we
indicate how to construct on any non-degenerate Legendre surface
$M\subset \mathcal{K}$ a canonical lift $\mathcal{A}:M\to
\widetilde{G}$ to the Lie sphere group (the \it{normal frame field
}\rm along $M$). By the means of the normal frame field we recover
the \it{Blaschke co-frame }\rm $(\alpha^1,\alpha^2)$ of $M$ and we
introduce a complete set of local differential invariants
$(q_1,q_2,p_1,p_2,r_1,r_2)$, the \it{invariant functions }\rm of
the surface. From the structural equations of the group
$\widetilde{G}$ we deduce the compatibility conditions fulfilled
by the normal co-frame and by the invariant functions. In $\S 3$
we analyze the \it{Lie-invariant Gauss map }\rm  and we write the
Euler-Lagrange equations of the variational problem in terms of
the invariant functions. Subsequently we set up a Pfaffian
differential system $(\mathcal{I},\Omega)$ on $P=G\times \R^6$
with the defining property that its integral manifolds are the
canonical lifts of Lie-minimal surfaces. In $\S 4$ we prove that
the differential system $(\mathcal{I},\Omega)$ is in involution
and that its general integral manifolds depend on six functions in
one variable. In the last part of the paper we prove our main
result :
\medskip

 $\mathbf{Theorem.}$ \it{Let $\Gamma\subset \mathcal{K}$
be a real-analytic Legendre curve and let $\mathcal{U}(\Gamma)\to
\Gamma$ be the corresponding tautological bundle. Let $L\subset
\mathcal{U}(\Gamma)$ be a real-analytic line sub-bundle of
$\mathcal{U}(\Gamma)$ and let $h,w:\Gamma\to \R$ be two
real-analytic functions. If $\Gamma$ and $L$ are suitably general,
then there exist a real-analytic Lie-minimal surface $M\subset
\mathcal{K}$ containing $\Gamma$ such that
$$\Gamma^*(\alpha^1+\alpha^2)=0,\quad
L=\Sigma_0(M)|_{\Gamma},\quad h=-3(q_1+q_2)|_{\Gamma},\quad
w=\frac{1}{3}(p_1-p_2)|_{\Gamma}.$$ Moreover, $M$ is unique in the
sense that any other Legendre surface with these properties agrees
with $M$ on an open neighborhood of $\Gamma$.}\rm
\medskip

 \section{Legendre Surfaces}\label{section1}

\subsection{Lie sphere geometry.} Let us begin with some basic
facts. Consider the vector space $\R^{(4,2)}$ with the inner
product of signature $(4,2)$ defined by
\begin{equation}
\langle V,V\rangle =
-2v^0v^5-2v^1v^4+(v^2)^2+(v^3)^2=g_{IJ}v^Iv^J\label{1.1.1}
\end{equation}
where $v^0,...,v^5$ are the components of $V$ with respect to the
standard basis $(\epsilon_0,...,\epsilon_5)$. We let $G$ be the
identity component of the pseudo-orthogonal group of (\ref{1.1.1})
and we let $\g$ be its Lie algebra. For each $A\in G$ we denote by
$A_J=A\cdot \epsilon_J$ the $J$-th column vector of $A$. Thus,
$(A_0,...,A_5)$ is a positive-oriented basis of $\R^{(4,2)}$ such
that $$\langle A_I,A_J\rangle =g_{IJ},\quad I,J=0,...,5.$$
Expressing the exterior derivative $dA_J$ in terms of the basis
$(A_0,...,A_5)$ we obtain
\begin{equation}
dA_J=\omega^I_JA_I,\quad J=0,...,5,\label{1.1.2}\end{equation}
where $\omega=(\omega^I_J)$ is the $\g$-valued Maurer-Cartan form
$A^{-1}dA$ on the group $G$. Taking the exterior derivative in
(\ref{1.1.2}) yields the \it{structure equations }\rm
\begin{equation}d\omega = -\omega\wedge \omega.\label{1.1.3}\end{equation}
The \it{Lie quadric }\rm is the space $\mathcal{Q}\subset \RP^5$
of the isotropic lines through $0\in \R^{(4,2)}$ and the \it{Lie
sphere group }\rm is defined to be the group
$\widetilde{G}=G/\{\pm I\}$ of all projective transformations of
$\RP^5$ which send $\mathcal{Q}$ into itself. Elements of
$\widetilde{G}$ are equivalence classes of matrices $A\in G$.
Given any such matrix $A$, its equivalence class in $G$ is denoted
by $[A]$. Thus, $[A]=[B]$ iff $A=\pm B$. Since the Maurer-Cartan
form $\omega=A^{-1}dA$ is bi-invariant under the action of $\{\pm
I\}$, then we can identify the Lie algebra of $\widetilde{G}$ with
$\g$ and we may think of $\omega$ as being the Maurer-Cartan form
of the Lie sphere group $\widetilde{G}$.

\begin{remark}\rm{The role of the Lie quadric is to represent the
oriented spheres of $\R^3$ (including point spheres, oriented
planes  and the "point at infinity"). Given a point $p\in \R^3$
and a real number $r$ we let $\sigma(p,r)$ be the oriented sphere
with center $p$ and signed radius $r$. Similarly, for every $p\in
\R^3$ and every $\overrightarrow{n}\in S^2$, we let
$\pi(p,\overrightarrow{n})$ be the oriented plane passing through
$p$ and orthogonal to the unit vector $\overrightarrow{n}$. Then,
the correspondence between the points of $\mathcal{Q}$ and
oriented spheres is given by
\begin{equation}
\left\{ \begin{array}{lll} \sigma(p,r)\to
\left[\left(1,\frac{r+p^1}{\sqrt{2}},p^2,p^3,\frac{r-p^1}{\sqrt{2}},\frac{p\cdot
p-r^2 }{2}\right)\right],\nonumber \\ \pi(p,\overrightarrow{n})\to
\left[\left(0,\frac{1+n^1}{2},\frac{n^2}{\sqrt{2}},\frac{n^3}{\sqrt{2}},\frac{1-n^1}{2},\frac{n\cdot
p }{\sqrt{2}}\right)\right],\\
 \infty \to [\epsilon_5]
\end{array}\right.\label{1.1.3.1}
\end{equation}}\end{remark}

\begin{defn}\rm{The \it{Kepler manifold }\rm $\mathcal{K}$ is defined to
be the isotropic Grassmannian consisting of all null planes
through the origin of $\R^{(4,2)}$.}\end{defn}

\begin{remark}\rm{Two oriented spheres corresponding to $[V],[V']\in
\mathcal{Q}$ are in oriented contact if and only if $\langle
V,V'\rangle =0$. Geometrically, this means that a null plane
$\ell\in \mathcal{K}$ represents a pencil of oriented spheres in
oriented contact. Thus, we may think of $\mathcal{K}$ as the
\it{space of the parabolic pencils of oriented spheres }\rm.
Another classical model of the Kepler manifold is "$\R^3\times
S^2$ with a $2$-dimensional sphere $S^2_{\infty}$ at the
infinity". The sphere $S^2_{\infty}$ is the set of the null planes
through the isotropic line $[\epsilon_5]$. The complement
$\mathcal{K}_0$ of $S^2_{\infty}$ is identified with $\R^3\times
S^2$ by
\begin{equation}F:(p,\overrightarrow{n})\in \R^3\times S^2\to [F_0(p)\wedge
F_1(p,\overrightarrow{n})]\in \mathcal{K},
\label{1.1.3.2}\end{equation} where
\begin{equation}
\left\{ \begin{array}{ll}
F_0(p)=\epsilon_0+\frac{p^1}{\sqrt{2}}\epsilon_1+p^2\epsilon_2+p^3\epsilon_3-\frac{p^1}{\sqrt{2}}\epsilon_4+
\frac{p\cdot p}{2}\epsilon_5,\\
F_1(p,\overrightarrow{n})=\frac{1+n^1}{2}\epsilon_1+\frac{n^2}{\sqrt{2}}\epsilon_2+\frac{n^3}{\sqrt{2}}\epsilon_3+\frac{1-n^1}{2}\epsilon_4+
\frac{n\cdot p }{\sqrt{2}}\epsilon_5.
\end{array}\right. \label{1.1.3.3}
\end{equation}}\end{remark}

Let $G$ act on $\mathcal{K}$ in the usual way : given a null plane
$[V\wedge V']$ spanned by a pair $V,V'$ of isotropic vectors, and
given $A\in G$, then $A\cdot [V\wedge V']=[AV\wedge AV']$. The
projection map
\begin{equation}\pi_{\mathcal{K}}:[A]\in G\to [A_0\wedge A_1]\in
\mathcal{K} \label{1.1.4}\end{equation} makes $G$ into a
$G_0$-principal fibre bundle over $\mathcal{K}$, where
$$G_0=\{[A]\in G : A^I_0=A^I_1=0,\quad I=2,...,5\}.$$
\begin{remark}\rm{For later
use we observe that the elements of $G_0$ can be written as
$$
X(D,B,Y,b)=\left(\begin{array}{ccc}
  D & D Y^*J  B & Z(D,Y,b) \\
  0 & B & Y \\
  0 & 0 & (D^{*})^{-1}
\end{array}\right),
$$
\noindent where $$D\in GL_{+}(2,\R),\quad B\in SO(2),\quad Y\in
\gl(2,\R),\quad b\in \R$$ \noindent and where
$$Z(D,Y,b)=\frac{1}{2}D J \left(Y^t Y + \left(\begin{array}{cc}
  0 & -b \\
  b & 0
\end{array}\right)\right),\quad J=\left(\begin{array}{cc}
  0 & 1 \\
  1 & 0
\end{array}\right),\quad X^*=J^tXJ.$$
}\end{remark}

The forms $\{\omega^2_0,
\omega^3_0,\omega^2_1,\omega^3_1,\omega^4_0\}$ are linearly
independent and span the semi-basic forms for the projection
$\pi_{\mathcal{K}}$. In particular, $\omega^4_0$ is well-defined
up to a positive multiple on $\mathcal{K}$ and, from the structure
equations, we get
$$d\omega^4_0=\omega^2_0\wedge\omega^2_1+\omega^3_0\wedge\omega^3_1+(\omega^0_0+\omega^1_1)\wedge
\omega^4_0.$$ From this we infer that $\omega^4_0$ defines a
$G$-invariant contact structure on $\mathcal{K}$.
\medskip

We let $\D$ be the \it{Dupin manifold }\rm, that is the manifold
of all $3$-dimensional linear subspaces of $\R^{(4,2)}$ of
signature $(2,1)$. Then, $G$ act transitively on $\D$ by $$A\cdot
[V\wedge V'\wedge V'']=[AV\wedge AV'\wedge AV''],$$ for every
$A\in G$ and every $[V\wedge V'\wedge V'']\in \D$. The projection
map $$\pi_{D}:A\in G\to [A_0\wedge A_3\wedge A_5]\in \D$$ gives on
$G$ the structure of a principal fibre bundle with structure group
$$S(\left( O(2,1)\times O(2,1) \right)\cong \{A\in G: A\cdot
[\epsilon_0\wedge \epsilon_3\wedge \epsilon_5]=[\epsilon_0\wedge
\epsilon_3\wedge \epsilon_5]\}.$$ Thus, $\D$ can be viewed as the
pseudo-riemannian symmetric space $$SO(4,2)/S\left(O(2,1)\times
O(2,1)\right).$$ The canonical pseudo-Riemannian metric $g_{D}$ on
$\D$ induced by the fibering $\pi_D$ is represented by the
tensorial quadratic form on $G$ given by
\begin{equation}
\omega^1_0\omega^0_1+\omega^0_4\omega^4_0+\omega^2_0\omega^0_2+2\omega^1_3\omega^3_1-\frac{1}{2}(\omega^3_2)^2.
\label{1.1.5}\end{equation}

\subsection{Legendre Surfaces}

\begin{defn}\rm{An oriented, connected immersed surface $M\subset \mathcal{K}$
is said to be \it{Legendrian }\rm if it is tangent to the contact
distribution on $\mathcal{K}$.} \end{defn}

Locally, there exist two smooth mappings $\phi_0,\phi_1:U\subset
M\to \R^{(4,2)}$ such that $\ell=[\phi_0(\ell)\wedge
\phi_1(\ell)]$, for every $\ell\in U$, and that $$\| \phi_0\|^2 =
\|\phi_1\|^2=\langle \phi_0,\phi_1\rangle = 0,\quad \langle
d\phi_0,\phi_1 \rangle = 0.$$ Then, $\langle
d\phi_0,d\phi_0\rangle$ and $\langle d\phi_1,d\phi_1 \rangle$ span
a sheaf $\mathcal{S}$ of quadratic forms on $M$. Throughout the
paper we shall assume that the fiber $\mathcal{S}_m$ of
$\mathcal{S}$ is two-dimensional, for every point $m\in M$. We
consider $$\mathcal{F}_0(M)=\{(\ell,A)\in M\times G \mid
\ell=[A_0\wedge A_1]\}$$ the pull back of the fiber bundle
$\pi_{\mathcal{K}}:G\to \mathcal{K}$ to the surface $M$. The local
cross sections of $\mathcal{F}_0(M)$ are called \it{local frame
fields }\rm along $M$. They can be considered as smooth maps
$A:U\to G$, where $U$ is an open subset of $M$, such that
$\ell=[A_0(\ell)\wedge A_1(\ell)]$, for every $\ell\in U$. For
every local frame field $A:U\to G$ we let $\alpha=(\alpha^I_J)$ be
the pull-back of the Maurer-Cartan form. Any other local frame
field $\widetilde{A}$ on $U$ is given by $\widetilde{A}=A\cdot X$,
where $X=U\to G_0$ is a smooth map. Thus, the $1$-forms $\alpha$
and $\widetilde{\alpha}$ are related by the gauge transformation
\begin{equation}
\widetilde{\alpha}=X^{-1}dX+X^{-1}\alpha X.
\label{1.2.1}\end{equation} A frame field $A:U\to G$ is of
\it{first order }\rm if, with respect to it
\begin{equation}\alpha^3_0\wedge \alpha^2_0>0,\quad
\alpha^2_0=\alpha^3_1=0.\label{1.2.2}\end{equation}From
(\ref{1.2.1}) it follows that first order frames exist on a
neighborhood of any point of $M$. The totality of first order
frames is a principal $G_1$-bundle $\mathcal{F}_1(M)\to M$ where
$$G_1=\left\{X(D,B,Y,b)\in G_1 : B=\epsilon Id_{2\times 2},\quad
D=\left(\begin{array}{cc}
 r & 0 \\
 0 & s
\end{array}\right),\quad \epsilon=\pm 1, rs>0 \right\}.$$

$\mathbf{Notation.}$  The elements of $G_1$ will be denoted by
$Y_{\epsilon}(r,s,Y,b)$, where $\epsilon = \pm1$, $Y\in
\gl(2,\R)$, $b,r,s\in \R$ and $rs>0$.
\medskip

We let $\sigma_0,\sigma_1:M\to \mathcal{Q}$ be defined by
$\sigma_0|_U=[A_0]$ and by $\sigma_1|_U=[A_1]$, for every first
order frame $A:U\to G$. We follow the classical terminology and we
call $\sigma_0$ and $\sigma_1$ the \it{curvature sphere maps }\rm
of the surface. We remark that $\sigma_0$ and $\sigma_1$ define a
splitting of the \it{tautological vector bundle }\rm
$$\mathcal{U}(M)=\{(\ell,V)\in M\times \R^{(4,2)}:\ell\in M,V\in
\ell\}$$ into the direct sum $\Sigma_0(M)\oplus \Sigma_1(M)$ of
the line sub-bundles $$\Sigma_a(M)=\{(\ell,V)\in M\times
\R^{(4,2)}:\ell\in M,V\in \sigma_a(\ell)\},\quad a=0,1.$$

\begin{defn}\rm{We say that $M$ is \it{non-degenerate }\rm if $\sigma_0$ and
$\sigma_1$ are immersions of $M$ into the Lie quadric.}\end{defn}

\begin{thm}{Let $M\subset \mathcal{K}$ be a non-degenerate
Legendre surface. Then there exist a unique lift $\mathcal{A}:M\to
\widetilde{G}$ to the group $\widetilde{G}$ satisfying the
Pfaffian equations
\begin{equation}\alpha^4_0=\alpha^2_0=\alpha^3_1=
\alpha^3_2=\alpha^1_0-\alpha^2_1=\alpha^0_1-\alpha^3_0=\alpha^0_2=\alpha^1_3=0,\label{PE1}
\end{equation}
\noindent and the independence condition
\begin{equation}\alpha^3_0\wedge \alpha^2_1>0.\label{PE2}
\end{equation} }\end{thm}
\begin{proof}\rm{If $A:U\to G$ is a first order frame
then, the linear differential forms $\alpha^1=\alpha^3_0$ and
$\alpha^2=\alpha^2_1$ give a positive-oriented co-framing on $U$
so that we  may write $$\alpha = P_1\alpha^1+P_2\alpha^2,$$ where
$P_1,P_2:U\to \g$ are smooth maps. The components of $P_a$ are
denoted by $P^I_{Ja}$, where $a=1,2$ and where $I,J=0,...,5$. If
$A,\widetilde{A}:U\to G$ are first order frames on $U$ and if the
corresponding transition function is of the form
$Y_{\epsilon}(r,s,Y,b):U\to G_1$, then
\begin{equation}
\tilde{\alpha}^1=\epsilon r\alpha^1,\quad
\tilde{\alpha}^2=\epsilon s\alpha^2,\quad
\tilde{\alpha}^1_0=r\left(s^{-1}\alpha^1_0-
Y^{2}_1\alpha^1\right),\quad
\tilde{\alpha}^0_1=s\left(r^{-1}\alpha^0_1-
Y^{1}_{2}\alpha^2\right).\label{FR1}
\end{equation}
This implies
\begin{equation}
\left\{ \begin{array}{ll} \widetilde{P}^{1}_{01}=\quad
\epsilon(s^{-1}P^{1}_{01}- Y^2_1),\quad
\widetilde{P}^{1}_{02}=\epsilon rs^{-2}P^{1}_{02},\\
\widetilde{P}^{0}_{12}=-\epsilon(r^{-1}P^{0}_{12}+Y^1_2),\quad
\widetilde{P}^{0}_{11}=\epsilon sr^{-2}P^{0}_{11}.
\end{array}\right.\label{FR2}
\end{equation}
From (\ref{FR2}) we see that for every point $\ell\in M$ there
exist a first order frame field $A:U\to G$ defined on an open
neighborhood $U$ of $\ell$ with respect to which
$P^1_{01}=P^0_{12}=0$. Such first order frame fields are said to
be of \it{second order }\rm. In addition, any other second order
frame field on $U$ is of the form $\widetilde{A}=A\cdot X$, where
$X:U\to G_2$ is a smooth map and
$$G_2=\left\{Y_{\epsilon}(r,s,Y,b)\in G_1 : Y
=\left(\begin{array}{cc}
  p & 0 \\
  0 & q
\end{array}\right),\quad p,q\in \R\right\}.$$

$\mathbf{Notation}$. The elements of $G_2$ will be denoted by
$Y_{\epsilon}(r,s,p,q,b)$, where $\epsilon = \pm 1$, $p,q,r,s,b\in
\R$ and $rs>0$.
\medskip

\noindent Second order frame fields are the cross sections of a
reduced sub-bundle $\mathcal{F}_2(M)$ of $\mathcal{F}_1(M)$ with
structural group $G_2$. Differentiating $\alpha^2_0=\alpha^3_1=0$
and applying the structure equations and Cartan's Lemma, we have
that $\alpha^3_2=0$, for every second order frame field $A$.
 Taking the exterior derivative of $\alpha^3_2=0$ and using again the structure equations
 and the Cartan's lemma we have $ \alpha^0_2\wedge \alpha^1-\alpha^1_3\wedge
 \alpha^2=0$. This implies $P^0_{22}=-P^1_{31}$ and hence we may write
\begin{equation}
\alpha^0_2=P^0_{21}\alpha^1+P^0_{22}\alpha^2,\quad
\alpha^1_3=-P^0_{22}\alpha^1+P^1_{32}\alpha^2.\label{FR3}
\end{equation}
If $A$ and $\widetilde{A}$ are second order frame fields on
$U\subset M$ and if $Y_{\epsilon}(r,s,p,q,b):U\to G_2$ is the
corresponding transition function, we then have \begin{equation}
\widetilde{\alpha}^0_2=\epsilon r^{-1}sp\alpha^0_1+\epsilon
r^{-1}\alpha^0_2-sb\alpha^2,\quad \widetilde{\alpha}^1_3=\epsilon
s^{-1}rq\alpha^1_0+\epsilon
s^{-1}\alpha^1_3+rb\alpha^1.\label{FR4}
\end{equation}
From this we obtain
\begin{equation}
\widetilde{P}^0_{21}=r^{-2}\left(P^{0}_{21}+spP^0_{11}\right),\quad
\widetilde{P}^{0}_{22}=\frac{1}{rs}P^0_{22}-\epsilon b, \quad
\widetilde{P}^1_{32}=s^{-2}\left(P^{1}_{32}+rqP^1_{02}\right).\label{FR5}
\end{equation}
Thus, for every point $\ell\in M$ there exist a second order frame
field $U\to G$ defined on an open neighborhood of $\ell$ with
respect to which
\begin{equation}P^{0}_{22}=P^{1}_{31}=0.\label{FR6}\end{equation} Such frame
fields are said to be of \it{third order }\rm. Now, (\ref{FR5})
implies that these frame fields are the local cross sections of a
reduced sub-bundle $\mathcal{F}_3(M)$ of $\mathcal{F}_2(M)$. The
structure group of $\mathcal{F}_3(\phi)$ is
$$G_3=\{Y_{\epsilon}(r,s,p,q,b)\in G_2 :b=0\}.$$ Since $M$ is
non-degenerate, then the functions $P^1_{02}$ and $P^0_{11}$ are
nowhere vanishing. Thus, from (\ref{FR2}) we infer that for every
$\ell\in M$ there exist a third order frame field $A:U\to G$
defined on an open neighborhood $U$ of $\ell$ such that
\begin{equation}\alpha^0_1=\alpha^1,\quad
\alpha^1_0=\alpha^2.\label{FR7}\end{equation} A third order frame
field satisfying (\ref{FR7}) is said to be of \it{fourth order
}\rm. If $A$ is a fourth order frame field on $U$, then any other
is given by $\widetilde{A}=AX$, where $X:U\to G_4$ and
$$G_4=\{Y_{\epsilon}(r,s,p,q,0)\in G_3 : r=s=\epsilon\}.$$ The
elements of $G_4$ are denoted by $Y_{\epsilon}(p,q)$, where
$p,q\in \R$. From this we immediately see that the third order
frame fields define a $G_4$ sub-bundle $\mathcal{F}_4(M)$ of
$\mathcal{F}_3(M)$. Now, using (\ref{FR5}) we see that for every
$\ell\in M$ there exist a fourth order frame $A:U\to G$ in an open
neighborhood $U$ of $\ell$ satisfying $$\alpha^0_2=\alpha^1_3=0.$$
Such frame fields are of \it{ fifth order }\rm. Notice that fifth
order frame satisfy(\ref{PE1}) and (\ref{PE2}).  Moreover, if $A$
and $\widetilde{A}$ are fifth order frame fields on an open
neighborhood $U$, then $\widetilde{A}=\epsilon A$, where $\epsilon
= \pm1$. This implies that the canonical lift $\mathcal{A}:M\to
\widetilde{G}$ is defined by $\mathcal{A}|_U=[A]$, for every fifth
order frame field $A:U\to G$.\qer }\end{proof}

\begin{defn}\rm{The map $\mathcal{A}:M\to \widetilde{G}$ is said to be the \it{normal frame field }\rm
along $M$. We then consider the linearly independent $1$-forms
$\alpha^1=\alpha^3_0$ and $\alpha^2=\alpha^2_1$ and we call
$(\alpha^1,\alpha^2)$ the \it{canonical co-frame }\rm of the
Legendre surface.}\end{defn}

\begin{remark}\rm{ If $M$ is the contact lift of $f:M\to \R^3$ and if $(u^1,u^2)$ are curvature lines coordinates
then $\alpha^1$ and $\alpha^2$ coincide with the Blaschke's
differentials
\begin{equation}
\left\{ \begin{array}{ll}
\alpha^1=\frac{1}{k_1-k_2}\left(\sqrt{\frac{g_{11}}{g_{22}}}(\partial_1k_1)^2\partial_2k_2)\right)^{\frac{1}{3}}du^1,\\
\alpha^2=\frac{-1}{k_1-k_2}\left(\sqrt{\frac{g_{22}}{g_{11}}}\partial_1k_1(\partial_2k_2)^2\right)^{\frac{1}{3}}du^2.
\end{array}\right.\end{equation}}\end{remark}

Taking the exterior derivatives of (\ref{PE1}) and using the
Maurer-Cartan equations it follows that there exist smooth
functions $q_1,q_2,p_1,p_2$ and $r_1,r_2$ such that
\begin{equation}
\left\{ \begin{array}{lll} \alpha^0_0=-2q_1\alpha^1+q_2\alpha^2, &
\alpha^1_1=-q_1\alpha^1+2q_2\alpha^2,\\
\alpha^0_3=r_1\alpha^1+p_2\alpha^2,&
\alpha^1_2=p_1\alpha^1+r_2\alpha^2,\\
\alpha^0_4=-r_2\alpha^1+r_1\alpha^2,
\end{array}\right. \label{1.2.6}
\end{equation}
We call $q_1,q_2,p_1,p_2$ and $r_1,r_2$ the \it{invariant
functions }\rm of the Legendre surface. Using once more the
Maurer-Cartan equations we obtain
\begin{equation}
d\alpha^1=\alpha^0_0\wedge\alpha^1,\quad
d\alpha^2=\alpha^1_1\wedge\alpha^2,\quad
d\alpha^1_2=-\alpha^1_1\wedge\alpha^1_2,\quad
d\alpha^0_3=-\alpha^0_0\wedge\alpha^0_3, \label{1.2.7}
\end{equation}
and
\begin{equation}
d\alpha^0_0=(\alpha^2-\alpha^0_3)\wedge\alpha^1,\quad
d\alpha^1_1=(\alpha^1-\alpha^1_2)\wedge\alpha^2,\quad
d\alpha^0_4=-(\alpha^0_0+\alpha^1_1)\wedge\alpha^0_4.
\label{1.2.8}
\end{equation}
We may rewrite these equations in terms of the invariant functions
\begin{equation}
d\alpha^1=-q_2\alpha^1\wedge \alpha^2,\quad
d\alpha^2=-q_1\alpha^1\wedge \alpha^2 \label{1.2.9}
\end{equation}
\begin{equation}
\left\{ \begin{array}{ll} -&2dq_1\wedge \alpha^1+dq_2\wedge
\alpha^2=(p_2-q_1q_2-1)\alpha^1\wedge\alpha^2,\\ -&dq_1\wedge
\alpha^1+2dq_2\wedge \alpha^2=(-p_1+q_1q_2+1)\alpha^1\wedge
\alpha^2,
\end{array}\right. \label{1.2.10}
\end{equation}
\begin{equation}
\left\{ \begin{array}{lll} & dr_1\wedge \alpha^1 + dp_2\wedge
\alpha^2 = (2q_2r_1+3q_1p_2)\alpha^1\wedge \alpha^2,\\ &
dp_1\wedge \alpha^1+dr_2\wedge \alpha^2 =
(2q_1r_2+3q_2p_1)\alpha^1\wedge \alpha 2,\\ -&dr_2\wedge
\alpha^1+dr_1\wedge \alpha^2 = 4(q_1r_1-q_2r_2)\alpha^1\wedge
\alpha^2.
\end{array}\right. \label{1.2.11}
\end{equation}

\section{Lie minimal surfaces}

\subsection{The Gauss map and the Euler-Lagrange equations}

\begin{defn}\rm{The \it{Gauss map }\rm of a generic Legendre surface $M\subset
\mathcal{K}$ is defined by $$\mathcal{D}:\ell\in M\to
[A_0(\ell)\wedge A_3(\ell)\wedge A_5(\ell)],\quad \forall \ell\in
M,$$ where $\mathcal{A}=[A]$ is the normal frame field along
$M$.}\end{defn}

\begin{remark}\rm{If we consider on $M$ the quadratic form
$\Phi=\alpha^1\alpha^2$, then (\ref{1.1.5}) and (\ref{PE1}) imply
that the Gauss map $\mathcal{D}:(M,\Phi)\to (\D,g_D)$ is an
isometric immersion.}\end{remark}

\begin{defn}\rm{A non-degenerate Legendre surface $M\subset
\mathcal{K}$ is said to be \it{Lie-minimal }\rm if it is a
critical point of the functional $$M\subset \mathcal{K}\to \int_M
\alpha^1\wedge \alpha^2$$ with respect to compactly supported
variations.}\end{defn}

\begin{remark}\rm{It is known (cfr. \cite{BH}) that Lie-minimal
surfaces are characterized by the harmonicity of the Gauss map.
}\end{remark}

\begin{thm}{Let $M\subset \mathcal{K}$ be a non-degenerate
Legendre surface. Then, $M$ is Lie-minimal if and only if
\begin{equation}
\left\{ \begin{array}{ll} dr_1\wedge
\alpha^2-4q_1r_1\alpha^1\wedge \alpha^2=0,\\ dr_2\wedge
\alpha^1-4q_2r_2\alpha^1\wedge \alpha^2=0.
\end{array}\right. \label{EL}
\end{equation}}\end{thm}

\begin{proof}\rm{ Without loss of generality we assume that
$\mathcal{D}$ is an embedding and we identify $M$ with its image
in $\D$. We extend the normal frame field $\mathcal{A}$ to a local
cross section $\widetilde{A}:\mathcal{U}\to G$ of $\pi_{D}:G\to
\D$ defined on an open neighborhood $\mathcal{U}\subset \D$ of
$M$. If we set
$\widetilde{\alpha}=\widetilde{\mathcal{A}}^{-1}d\widetilde{\mathcal{A}}$
then
\begin{equation}
\left\{ \begin{array}{ll} \beta^1=\widetilde{\alpha}^1_0,\quad
\beta^2=\widetilde{\alpha}^0_1,\quad
\beta^3=\widetilde{\alpha}^0_4,\quad
\beta^4=\widetilde{\alpha}^4_0,\quad
\beta^5=\widetilde{\alpha}^2_0\\
\beta^6=\widetilde{\alpha}^0_2,\quad
\beta^7=\widetilde{\alpha}^1_3,\quad
\beta^8=\widetilde{\alpha}^3_1,\quad
\beta^9=\widetilde{\alpha}^3_2,
\end{array}\right. \label{2.1.1}
\end{equation}
is a co-frame on $\mathcal{U}$. We let $B_1,...,B_9$ be the local
trivialization of $T(\D)$ dual to $(\beta^1,...,\beta^9)$. Then,
(\ref{1.2.6}) implies that
\begin{equation}X_1=B_2|_M-r_2B_3|_M,\quad
X_2=B_1|_M+r_1B_3|_M\label{2.1.2}\end{equation} is the
trivialization\footnote{Here $T(M)$ is viewed as a sub-bundle of
$\mathcal{D}^*\left(T(\D)\right)$} of $T(M)$ dual to the canonical
co-frame $(\alpha^1,\alpha^2)$ and that
\begin{equation}\overline{B}_3=B_3|_M\quad \overline{B}_4=
B_4|_M+r_2B_1|_M-r_1B_2|_M,\quad
\overline{B}_5=B_5|_M,....,\overline{B}_9=B_9|_M\label{2.1.3}\end{equation}
is a trivialization of the normal bundle $\mathcal{N}\to M$.
Taking the exterior derivative of (\ref{2.1.1}) and using the
Maurer-Cartan equations we compute the covariant derivatives of
the vector fields $B_1,...,B_9$. We then have
\begin{equation}
\left\{ \begin{array}{lll} \nabla
B_1=(\widetilde{\alpha}^1_1-\widetilde{\alpha}^0_0)B_1+\widetilde{\alpha}^2_1B_5-\widetilde{\alpha}^0_3B_7,\\
\nabla
B_2=(\widetilde{\alpha}^0_0-\widetilde{\alpha}^1_1)B_2-\widetilde{\alpha}^1_2B_6+\widetilde{\alpha}^3_0B_8,\\
\nabla
B_3=(\widetilde{\alpha}^0_0+\widetilde{\alpha}^1_1)B_3-\widetilde{\alpha}^2_1B_6+\widetilde{\alpha}^3_0B_7,
\end{array}\right. \label{2.1.4}
\end{equation}
and
\begin{equation}
\left\{ \begin{array}{llllll} \nabla
B_4=-(\widetilde{\alpha}^0_0+\widetilde{\alpha}^1_1)B_4+\widetilde{\alpha}^1_2B_5-\widetilde{\alpha}^0_3B_8,\\
\nabla B_5=\widetilde{\alpha}^1_2B_1+\widetilde{\alpha}^2_1B_4
-\widetilde{\alpha}^0_0B_5+\widetilde{\alpha}^0_3B_9,\\ \nabla
B_6=-\widetilde{\alpha}^2_1B_2-\widetilde{\alpha}^1_2B_3
+\widetilde{\alpha}^0_0B_6+\widetilde{\alpha}^3_0B_9,\\ \nabla
B_7=-\widetilde{\alpha}^3_0B_1+\widetilde{\alpha}^0_3B_3
+\widetilde{\alpha}^1_1B_7-\widetilde{\alpha}^2_1B_9,\\ \nabla
B_8=\widetilde{\alpha}^0_3B_2-\widetilde{\alpha}^3_0B_4
-\widetilde{\alpha}^1_1B_8-\widetilde{\alpha}^1_2B_9,\\ \nabla
B_9=\widetilde{\alpha}^3_0B_5+\widetilde{\alpha}^0_3B_6
-\widetilde{\alpha}^1_2B_8-\widetilde{\alpha}^2_1B_9.\\
\end{array}\right. \label{2.1.5}
\end{equation}
Thus, from (\ref{2.1.2}), (\ref{2.1.3}), (\ref{2.1.4}) and
(\ref{2.1.5}) we infer that the shape operator $$S\in \Gamma(M,
\mathrm{Hom}(TM,\Omega^1(M)\otimes \mathcal{N}))$$ of $M\subset
\D$ is given by
\begin{equation}
\left\{ \begin{array}{ll}
S(X_1)=-\left(dr_2+2r_2(-q_1\alpha^1+2q_2\alpha^2)\right)\overline{B}_3+
\alpha^1\left(-p_1\overline{B}_6-r_2\overline{B}_7+\overline{B}_8\right),\\
S(X_2)=\left(dr_1+2r_1(-2q_1\alpha^1+q_2\alpha^2)\right)\overline{B}_3+\alpha^2\left(\overline{B}_5-r_1\overline{B}_6-
p_2\overline{B}_7\right).
\end{array}\right. \label{2.1.6}
\end{equation}
In particular, we obtain the following formula for the mean
curvature vector
\begin{equation}
H=\frac{1}{2}\left(S(X_1)(X_2)+S(X_2)(X_1)\right)=\left(-dr_2(X_2)-4r_2q_2+dr_1(X_1)-4p_1q_1\right)\overline{B}_3
\label{2.1.7}\end{equation} From this we deduce that
$\mathcal{D}:(M,\Phi)\to (\D,g_D)$ is harmonic if and only if
\begin{equation}
-dr_2\wedge \alpha^1+4q_2r_2\alpha^1\wedge \alpha^2+dr_1\wedge
\alpha^2-4r_1q_1\alpha^1\wedge\alpha^2=0.\label{2.1.8}\end{equation}
From (\ref{1.2.11}) and (\ref{2.1.8}) we get the required
result.\qer }\end{proof}

\subsection{The differential system of Lie-minimal surfaces}
Let $P$ be the \it{configuration space }\rm $\widetilde{G}\times
\R^6$ and let denote by $(q_1,q_2,p_1,p_2,r_1,r_2)$ the
coordinates on $\R^6$. On $P$ we consider the Pfaffian ideal
$\mathcal{I}\subset \Omega^{*}(P)$ generated (as a differential
ideal) by the linear differential forms
\begin{equation}
\left\{ \begin{array}{ll} & \eta^1=\omega^4_0,\quad
\eta^2=\omega^2_0,\quad \eta^3=\omega^3_1,\\ &
\eta^4=\omega^3_2,\quad \eta^5=\omega^1_0-\omega^2,\quad
\eta^6=\omega^0_1-\omega^1,\\ &\eta^7=\omega^0_2,\quad
\eta^8=\omega^1_3,
\end{array}\right. \label{2.2.1}
\end{equation}
\begin{equation}
\left\{ \begin{array}{ll} &
\eta^9=\omega^0_0+2q_1\omega^1-q_2\omega^2,\\ &
\eta^{10}=\omega^1_1+q_1\omega^1-2q_2\omega^2,
\end{array}\right. \label{2.2.2}
\end{equation}
\begin{equation}
\left\{ \begin{array}{lll} &
\eta^{11}=\omega^0_3-r_1\omega^1-p_2\omega^2,\\ &
\eta^{12}=\omega^1_2-p_1\omega^1-r_2\omega^2,\\ &
\eta^{13}=\omega^0_4+r_2\omega^1-r_1\omega^2,
\end{array}\right. \label{2.2.3}
\end{equation}
and by the exterior differential $2$-forms
\begin{equation}
\Theta^1=dr_1\wedge \alpha^2-4q_1r_1\alpha^1\wedge \alpha^2,\quad
\Theta^2=dr_2\wedge \alpha^1-4q_2r_2\alpha^1\wedge \alpha^2,
\label{2.2.4}
\end{equation}
together with the independence condition $\Omega=\omega^1\wedge
\omega^2$, where $\omega^1=\omega^3_0$ and $\omega^2=\omega^2_1$.
Taking the exterior derivatives of
(\ref{2.2.1}),(\ref{2.2.2}),(\ref{2.2.3}) and using the structural
equations of $\widetilde{G}$ we obtain the \it{quadratic equations
}\rm  \begin{equation} d\eta^1\equiv .....\equiv d\eta^8 \equiv
d\eta^{13}\equiv 0,\label{2.2.5}\end{equation} and
\begin{equation}
\left\{ \begin{array}{lll} & d\eta^9\equiv
2\pi^1\wedge\omega^1-\pi^2\wedge\omega^2+(-1+p_2-q_1q_2)\omega^1\wedge\omega^2,\\
&d\eta^{10}\equiv\pi^1\wedge\omega^1-2\pi^2\wedge\omega^2+(1-p_1+q_1q_2)\omega^1\wedge\omega^2,\\
&
d\eta^{11}\equiv\zeta^1\wedge\omega^1+\upsilon^2\wedge\omega^2-(2r_1q_2+3q_1p_2)\omega^1\wedge
\omega^2,
\\ & d\eta^{12}\equiv
\upsilon^1\wedge\omega^1+\zeta^2\wedge\omega^2-(3p_1q_2+2r_2q_1)\omega^1\wedge
\omega^2,
\end{array}\right.
\label{2.2.6}
\end{equation}
where $\equiv$ denotes equality up to the algebraic ideal
generated by $\eta^1,....,\eta^{13},\Theta^1,\Theta^2$ and where
$$\pi^i=dq_i,\quad \upsilon^i=dp_i,\quad \zeta^i=dr_i,\quad
i=1,2.$$ If we set
\begin{equation}
\left\{ \begin{array}{lll} &
\Omega^1=2\pi^1\wedge\omega^1-\pi^2\wedge\omega^2+(-1+p_2-q_1q_2)\omega^1\wedge\omega^2,\\
&\Omega^2=\pi^1\wedge\omega^1-2\pi^2\wedge\omega^2+(1-p_1+q_1q_2)\omega^1\wedge\omega^2,\\
&
\Omega^3=\zeta^1\wedge\omega^1+\upsilon^2\wedge\omega^2-(2r_1q_2+3q_1p_2)\omega^1\wedge
\omega^2,\\
 & \Omega^4=
\upsilon^1\wedge\omega^1+\zeta^2\wedge\omega^2-(3p_1q_2+2r_2q_1)\omega^1\wedge
\omega^2.
\end{array}\right.
\label{2.2.7}
\end{equation}
then
\begin{equation}\{\eta^1,....,\eta^{13},\Theta^1,\Theta^2,\Omega^1,\Omega^2,\Omega^3,\Omega^4,d\Theta^1,d\Theta^2\}
\label{2.2.8}
\end{equation}
is a set of algebraic generators of the differential ideal
$\mathcal{I}$. The integral manifolds of this system are
two-dimensional submanifolds $\widetilde{M} \subset P$ such that
$$ \eta^a=0,\quad
\Theta^1=\Theta^2=\Omega^1=\Omega^2=\Omega^3=0=\Omega^4,\quad
\Omega\neq 0.$$ Thus, the map
\begin{equation}\phi:([A],q_1,q_2,p_1,p_2,r_1,r_2)\in \widehat{M}\to [A_0\wedge A_1]\in
\mathcal{K}\label{2.2.9}\end{equation} is a non-degenerate
Legendre immersion. Since our arguments are local, we identify
$\widetilde{M}$ with its image $M=\phi(\widetilde{M})\subset
\mathcal{K}$. Thus, $M$ a Lie-minimal surface with normal frame
field
$$\mathcal{A}:([A],q_1,q_2,p_1,p_2,r_1,r_2)\in M\to [A]\in
\widetilde{G}.$$ Conversely, if $M$ is a Lie-minimal surface with
normal frame field $\mathcal{A}:M\to G$ and with invariant
functions $q_1,...,r_2$, then the map
\begin{equation}\ell\in M\to
\left(\mathcal{A}(\ell),q_1(\ell),q_2(\ell),p_1(\ell),p_2(\ell),r_1(\ell),r_2(\ell)\right)\in
P,\quad \forall \ell\in M\label{2.2.10}\end{equation} defines an
integral manifold of the differential system
$(\mathcal{I},\Omega)$. To summarize :

\begin{prop}{Lie-minimal surfaces $M\subset \mathcal{K}$ may be
regarded as being the integral submanifolds of the differential
system $(\mathcal{I},\Omega)$ on $P$.}\end{prop}

\section{The Cauchy problem}
\subsection{Involutivity of the differential system} On $P$ we consider the parallelization
$$\left(\frac{\partial}{\partial
\omega^i},\frac{\partial}{\partial
\eta^a},\frac{\partial}{\partial \pi^i},\frac{\partial}{\partial
\upsilon^i},\frac{\partial}{\partial \zeta}\right),\quad
i=1,2,a=1,...,13$$ dual to the co-frame
$(\omega^i,\eta^a,\pi^i,\upsilon^i,\zeta^i)$. We define
\begin{equation}
\left\{ \begin{array}{ll} V_1(\Omega)=\{(\z,E_1)\in G_1(T(P)):
\left((\omega^1)^2+(\omega^2)^2\right)|_{E_{1}}\neq 0\},\\
V_2(\Omega)=\{(\z,E_2)\in G_2(T(P)): \Omega|_{E_{2}}\neq 0\}
\end{array}\right. \label{3.1.1}
\end{equation}
and we set $$\V_1=S^1\times \R^{16},\quad \V_2=\R(13,2)\oplus
\R(2,2)\oplus \R(2,2)\oplus \R(2,2)\cong \R^{38},$$ with
coordinates $$z=(\cos(\theta),\sin(\theta),x^a,y^i,u^i,v^i),\quad
Z=(X^a_j,Y^i_j,U^i_j,V^i_j),\quad a=1,...,13, i,j=1,2.$$ Then, we
identify $V_1(\Omega)$ with $P\times \V_1$ and $V_2(\Omega)$ with
$P\times \V_2$ by the means of
\begin{equation}
\left\{ \begin{array}{ll} (\z,z)\in P\times \V_1\to
(\z,E_1(\z,z))\in V_1(\Omega),\\ (\z,Z)\in P\times V_2\to
(\z,E_2(\z,Z))\in V_2(\Omega),
\end{array}\right.
\label{3.1.2}
\end{equation}
where $$E_1(\z,z)=\left[\cos(\theta)\frac{\partial}{\partial
\omega^1}|_{\z}+ \sin(\theta)\frac{\partial}{\partial
\omega^2}|_{\z} +x^a\frac{\partial}{\partial \eta^a}|_{\z}+y^i
\frac{\partial}{\partial \pi^i}|_{\z}+u^i\frac{\partial}{\partial
\upsilon^i}|_{\z}+v^i\frac{\partial}{\partial
\zeta^i}|_{\z}\right]$$ and where $$ \left\{ \begin{array}{lll}
E_2(\z,Z)=\left[T_1(\z,Z)\wedge T_2(\z,Z)\right],\\
T_1(\z,Z)=\frac{\partial}{\partial \omega^1}|_{\z}
+X^a_1\frac{\partial}{\partial \eta^a}|_{\z}+Y^i_1
\frac{\partial}{\partial
\pi^i}|_{\z}+U^i_1\frac{\partial}{\partial
\upsilon^i}|_{\z}+V^i_1\frac{\partial}{\partial \zeta^i}|_{\z},\\
T_2(\z,Z)=\frac{\partial}{\partial \omega^2}|_{\z}
+X^a_2\frac{\partial}{\partial \eta^a}|_{\z}+Y^i_2
\frac{\partial}{\partial
\pi^i}|_{\z}+U^i_2\frac{\partial}{\partial
\upsilon^i}|_{\z}+V^i_2\frac{\partial}{\partial \zeta^i}|_{\z}.
\end{array}\right.
$$ Thus, the space $V_1(\mathcal{I},\Omega)$ consisting of the
$1$-dimensional integral elements can be identified with the
submanifold of $P\times \V_1$ defined by the linear equations
$x^1=....=x^{13}=0$. Similarly, the space
$V_2(\mathcal{I},\Omega)$ consisting of all $2$-dimensional
integral elements is identified with the submanifold of $P\times
\V_2$ defined by  $$ \left\{
\begin{array}{ll} X^a_i=0,\quad a=1,...,13, i=1,2,\\
2Y^1_2+Y^2_1+(1-p_2+q_1q_2)=Y^1_2+2Y^2_1-(1-p_1+q_1q_2)=0,\\
U^2_1-V^1_2-(2r_1q_2+3q_1p_2)=U^1_2-V^2_1+(3p_1q_2+2r_2q_1)=0,\\
V^1_1-4q_1r_1=V^2_2-4q_2r_2=0.
\end{array}\right.
$$ From this we infer that all the integral elements of the system
are $K$-ordinary. Let $(\z,E_1)$ be a $1$-dimensional integral
element such that $$E_1=\left[a^1\frac{\partial}{\partial
\omega^1}|_{\z}+ a^2\frac{\partial}{\partial \omega^2}|_{\z}+y^i
\frac{\partial}{\partial \pi^i}|_{\z}+u^i\frac{\partial}{\partial
\upsilon^i}|_{\z}+v^i\frac{\partial}{\partial
\zeta^i}|_{\z}\right].$$ Contracting $\Theta^1,\Theta^2$ and
$\Omega^1,...,\Omega^4$ with $E_1$, we obtain the polar equations
of the integral element $E_1$
\begin{equation}
\left\{ \begin{array}{lllllll} \eta^a=0,\quad a=1,...,13,\\
2a^1\pi^1-a^2\pi^2-\left(2y^1+a^2(1-p_2+q_1q_2)\right)\omega^1+\left(y^2+a^1(1-p_2+q_1q_2)\right)\omega^2=0,\\
a^1\pi^1-2a^2\pi^2-\left(y^1-a^2(1-p_1+q_1q_2)\right)\omega^1+\left(2y^2-a^1(1-p_1+q_1q_2)\right)\omega^2=0,\\
a^1\upsilon^1+a^2\zeta^2-\left(u^1+a^2(3p_1q_2+2r_2q_1)\right)\omega^1-\left(v^2-a^1(3p_1q_2+2r_2q_1)\right)\omega^2=0,\\
a^2\upsilon^2+a^1\zeta^1-\left(v^1+a^2(2r_1q_2+3q_1p_2)\right)\omega^1-\left(u^2-a^1(2r_1q_2+3q_1p_2)\right)\omega^2=0,\\
a^2\zeta^1-4a^2q_1r_1\omega^1+(4a^1q_1r_1+v^1)\omega^2=0,\\
a^1\zeta^2-(4a^2q_2r_2+v^2)\omega^1+4a^1q_2r_2\omega^2=0.
\end{array}\right.
\label{3.1.3}
\end{equation}
This shows that the polar space $H(\z,E_1)$ of a
\it{non-characteristic integral element}\rm\footnote{i.e. an
integral element such that $a^1a^2\neq 0$} is two-dimensional.
From the Cartan-Kaehler theorem we deduce

\begin{prop}{If $\widehat{\Gamma} \subset P$ is a non-characteristic real-analytic integral
curve of $(\mathcal{I},\Omega)$ then there exist a unique
real-analytic integral manifold $\widehat{M}\subset P$ such that
$\widehat{\Gamma}\subset \widehat{M}$.}\end{prop}

\subsection{Legendre Curves.}\label{Legendre Curves}
Let $\Gamma\subset \mathcal{K}$ be a smooth Legendre curve.
Locally,  we have that $\ell=[V_0(\ell)\wedge V_1(\ell)]$, for
every $\ell\in \Gamma$, where $V_0,V_1:\Gamma\to \R^{(4,2)}$ are
smooth maps such that
\begin{equation}\parallel V_0,\parallel =\parallel V_1\parallel =\langle
V_0,V_1\rangle=0,\quad \langle
V_0,dV_1\rangle=0.\label{LC1}\end{equation} We say that $\Gamma$
is \it{linearly full }\rm in case\footnote{ We use the notation
$dV=V'd\zeta$, where $d\zeta$ is a nowhere vanishing $1$-form on
$N$.}
\begin{equation}
V_0(\ell)\wedge V_1(\ell)\wedge V'_0(\ell)\wedge V'_1(\ell)\wedge
V''_0(\ell)\wedge V''_1(\ell)\neq 0\quad \forall \ell\in
N.\label{LC2}\end{equation} We let $\mathcal{U}(\Gamma)\to \Gamma$
be the tautological vector bundle of the curve, that is
\begin{equation}
\mathcal{U}(\Gamma)=\{(\ell,V)\in \Gamma\times \R^{(4,2)}:V\in
\ell\}.\label{LC3}\end{equation} A cross section of
$\mathcal{U}(\Gamma)$ is a smooth map $V:\Gamma\to \R^{(4,2)}$
such that $V(\ell)\in \ell$, for every $\ell\in \Gamma$.
Accordingly, a line sub-bundle $L\subset \mathcal{U}(\Gamma)$ can
be viewed as a mapping $\sigma_L:\Gamma \to \mathcal{Q}$ such that
$\sigma_L(\ell)\subset \ell$, for each $\ell\in \Gamma$. We say
that $L\subset \mathcal{U}(\Gamma)$ is \it{fat }\rm if
\begin{equation}
V(\ell)\wedge V'(\ell)\wedge .... \wedge V^{(v)}(\ell)\neq 0,\quad
\forall \ell\in \Gamma,\label{LC4}\end{equation} for every local
trivialization $V:U\to \R^{(4,2)}$ of $L$. In this case the
osculating space $\delta_L(\ell)=[V(\ell)\wedge V'(\ell)\wedge
V''(\ell)]\subset \R^{(4,2)}$ has signature $(2,1)$, for every
$\ell\in \Gamma$. The map
\begin{equation}
\delta_L:\ell\in \Gamma\to \delta_L(\ell)\in \D,\quad \forall
\ell\in \Gamma. \label{LC5}\end{equation} is called the
\it{directrix curve }\rm of $L$. If $\delta_L$ is non-isotropic
(i.e. $\delta_L^*(g_D)$ is nowhere vanishing) then $L$ is said to
be a \it{polarization }\rm of the curve $\Gamma$.

\begin{prop}{Let $(\Gamma,L)$ be a polarized Legendre curve. Then, there
exist a unique map $\mathcal{R}:\Gamma \to \widetilde{G}$ such
that
\begin{equation}\ell=[\mathrm{R}_0(\ell)\wedge\mathrm{R}_1(\ell)],\quad
\mathrm{R}_0(\ell)\in L|_{\ell},\quad \forall \ell\in
\Gamma\label{LC6}\end{equation} and that
\begin{eqnarray}
 \mathcal{R}^{-1}d\mathcal{R}&=&\left(\begin{array}{cccccc}
         k_0&1&0&k_1&k_3&0\\-1&-k_0&k_2&0&0&-k_3\\
         0&-1&0&0&k_2&0\\1&0&0&0&0&k_1\\0&0&-1&0&k_0&-1\\0&0&0&1&1&-k_0
\end{array}\right)\mu,\label{LC7}
\end{eqnarray}
where $\mu$ is a nowhere vanishing $1$-form and where
$k_0,k_1,k_2$ and $k_3$ are real-valued functions.}\end{prop}
\begin{proof}\rm{ We consider the $G_0$ fiber bundle
$$\mathcal{R}_0(\Gamma,L)=\{(\ell,\mathrm{R})\in \Gamma\times G :
\ell=[\mathrm{R}_0\wedge\mathrm{R}_0],\quad \mathrm{R}_0\in
L|_{\ell}\}.$$ The cross-sections of $\mathcal{R}_0(\Gamma,L)$ are
smooth maps $R:U \to G$ defined on an open subset $U\subset
\Gamma$, such that
$$\ell=[\mathrm{R}_0(\ell)\wedge\mathrm{R}_1(\ell)],\quad
\mathrm{R}_0(\ell)\in L|_{\ell},\quad \forall \ell\in U.$$ For
each frame field $R:U\to G$ we let $\rho$ be the $\g$-valued
$1$-form $\mathrm{R}^{-1}d\mathrm{R}$. We say that $R:U\to G$ is
of \it{first order }\rm if
 \begin{equation}\rho^3_0\neq
0,\quad
\rho^2_0=\rho^3_1=\rho^3_0+\rho^2_1=\rho^4_0=0.\label{LC8}\end{equation}
Since $\Gamma$ is linearly full then, first order frames do exist
near any point of $\Gamma$ and they define a sub-bundle
$\mathcal{R}_1(\Gamma,L)$ of $\mathcal{R}_0(\Gamma,L)$ with fiber
$$H_1=\{X\in G_0 : X =X(r\epsilon I,\epsilon I,Y,b):
\epsilon=\pm1, r,b\in \R, r\neq 0, Y\in \gl(2,\R) \}.$$ If
$\mathrm{R}$ and $\widetilde{\mathrm{R}}$ are first order frames
such that $\widetilde{\mathrm{R}}=\mathrm{R}X(r\epsilon I,\epsilon
I,Y,b)$ then
\begin{equation}\widetilde{\rho}^0_1=\rho^0_1+rY^2_1\rho^3_0,\quad
\widetilde{\rho}^1_0=\rho^1_0-rY^1_2\rho^3_0,\quad
\widetilde{\rho}^3_2=\rho^3_2+\epsilon
r(Y^1_2+Y^2_1)\rho^3_0.\label{LC9}\end{equation} This shows that
near any point of $\Gamma$ there exist first order frames such
that
\begin{equation}\rho^1_0+\rho^0_1=\rho^3_2=0.\label{LC10}\end{equation}
Frame fields satisfying (\ref{LC10}) are said to be of \it{second
order }\rm. From (\ref{LC9}) it follows that the totality of
second order frames defines a fiber bundle
$\mathcal{R}_2(\Gamma,L)$ with fiber $$H_2=\{X =X(r\epsilon
I,\epsilon I,Y,b)\in H_1: Y^1_2=Y^2_1=0 \}.$$ Notice that the
$1$-form $\rho^0_1$ is independent on the choice of the second
order frame and hence there exist $\mu\in \Omega^1(\Gamma)$ such
that $\mu|_U=\rho^0_1$. At this juncture it is convenient to
recall that the pseudo-riemannian metric $g_{D}$ of $\D$ is
represented by the tensorial quadratic form on $G$ defined by $$
\omega^1_0\omega^0_1+\omega^0_4\omega^4_0+\omega^2_0\omega^0_2+2\omega^1_3\omega^3_1-\frac{1}{2}(\omega^3_2)^2.
$$ Thus, $\delta_L^*(g_{D})=-\mu^2$ and hence $\mu$ is nowhere
vanishing. From (\ref{LC9}) it follows that, locally, there exist
second order frames such that
\begin{equation}
\rho^3_0=-\rho^2_1=-\rho^1_0=\rho^0_1=\mu.\label{LC11}\end{equation}
Frame fields satisfying (\ref{LC11}) are of \it{third order }\rm.
The totality of third order frames originates a principal fiber
bundle $\mathcal{R}_3(\Gamma,L)$ with structural group $$H_3=\{X
=X(r\epsilon I,\epsilon I,Y,b)\in H_2: r=1 \}.$$ If
$\widetilde{\mathrm{R}}$ and $\mathrm{R}$ are third order frame
fields then $$\widetilde{\rho}^0_0=\rho^0_0-\epsilon
Y^2_2\rho^3_0,\quad \widetilde{\rho}^1_1=\rho^1_1+\epsilon Y^1_1
\rho^3_0.$$ Therefore, near any point of $\Gamma$ there exist a
third order frame field $R$ such that
\begin{equation}\rho^1_1+\rho^0_0=0.\label{LC12}\end{equation} Frame fields satisfying (\ref{LC12}) define
a reduced sub-bundle $\mathcal{R}_4(\Gamma,L)$ with structure
group $$H_4=\{X =X(\epsilon I,\epsilon I,Y,b)\in H_3:
Y^1_1=Y^2_2\}.$$ Consider two local cross sections $R$ and
$\widetilde{R}$ of $\mathcal{R}_4(\Gamma,L)$, we then have
\begin{equation}\widetilde{\rho}^0_2=\rho^0_2+\epsilon
(Y^1_1+\frac{b}{2})\mu,\quad
\widetilde{\rho}^1_3=\rho^1_3-\epsilon
(Y^1_1-\frac{b}{2})\mu.\label{LC13}\end{equation} This implies
that there exist fourth order frame fields with respect to which
\begin{equation}\rho^0_2=\rho^1_3=0.\label{LC14}\end{equation}
Fourth order frame fields satisfying (\ref{LC14}) are said to be
of \it{fifth order }\rm. The totality of fifth order frames
generates a reduced sub-bundle $\mathcal{R}_5(\Gamma,L)$ with
fiber $\Z_2=\{\pm I\}$ and henceforth there exist a map
$\mathcal{R}:\Gamma\to \widetilde{G}$ such that
$\mathcal{R}|_U=[R]$, for every fifth order frame field $R:U\to
G$. From (\ref{LC8}), (\ref{LC10}),(\ref{LC11}),(\ref{LC12}) and
(\ref{LC14}) it follows that $\mathcal{R}$ satisfies the required
properties.\qer }\end{proof}

\begin{defn}\rm{The lift $\mathcal{R}:N\to
\widetilde{G}$ is said to be the \it{Frenet frame field }\rm of
$(\Gamma,L)$. The $1$-form $\mu$ is the \it{Lie-invariant line
element }\rm and the functions $k_0,k_1,k_2$ and $k_3$ are the
\it{generalized curvatures }\rm of $(\Gamma,L)$.}\end{defn}

\begin{remark}\rm{This proposition shows that polarized Legendre
curves are completely determined, up to the action of the Lie
sphere group, by the curvatures $k_0,...,k_3$.}\end{remark}

\subsection{The Cauchy problem}
\begin{thm}{Let $(\Gamma,L)$ be a real-analytic
polarized Legendre curve and let $h,w:\Gamma\to \R$ be two
real-analytic functions. Then, there exist a real-analytic
Lie-minimal surface $M\subset \mathcal{K}$ containing $\Gamma$
such that $$\Gamma^*(\alpha^1+\alpha^2)=0,\quad
L=\Sigma_0(M)|_{\Gamma},\quad h=-3(q_1+q_2)|_{\Gamma},\quad
w=\frac{1}{3}(p_1-p_2)|_{\Gamma}.$$ This manifold is unique in the
sense that any other Legendre surface with these properties agrees
with M on an open neighborhood of
$\Gamma$.}\label{maintheorem}\end{thm}
 \begin{proof}\rm{Let $\mathcal{R}:\Gamma \to \widetilde{G}$ be
 the Frenet frame field along $(\Gamma,L)$. We set
\begin{eqnarray}
 X(h)&=&\left(\begin{array}{cccccc}
         1&0&0&-h/2&h/2&h^2/8\\
         0&1&h/2&0&h^2/8&-h/2\\
         0&0&1&0&h/2&0\\
         0&0&0&1&0&h/2\\
         0&0&0&0&1&0\\
         0&0&0&0&0&1
\end{array}\right).\label{CP1}
\end{eqnarray}
and we consider the frame field
\begin{equation}
\widehat{\mathcal{R}}=\mathcal{R}X(h):\Gamma\to \widetilde{G}
.\label{CP2}\end{equation} We then have
\begin{eqnarray}
\widehat{\mathcal{R}}^{-1}d\widehat{\mathcal{R}}&=&\left(\begin{array}{cccccc}
         k_0+h/2&1&0&\widehat{k}_1&\widehat{k}_3&0\\-1&-k_0+h/2&\widehat{k}_2&0&0&-\widehat{k}_3\\
         0&-1&0&0&\widehat{k}_2&0\\1&0&0&0&0&\widehat{k}_1\\0&0&-1&0&k_0+h/2&-1\\0&0&0&1&1&-k_0+h/2
\end{array}\right)\mu,\label{CP3}
\end{eqnarray}
where \begin{equation} \left\{ \begin{array}{lll}
\widehat{k}_2=k_2+\frac{1}{2}h'-\frac{1}{2}(hk_0-\frac{1}{4}h^2),\\
\widehat{k}_1=k_1-\frac{1}{2}h'-\frac{1}{2}(hk_0+\frac{1}{4}h^2),\\
\widehat{k}_3=k_3-\frac{1}{2}h'-\frac{1}{4}h^2.
\end{array}\right.\label{CP4}
\end{equation}
We define $\overline{q}_i,\overline{p}_i,\overline{r}_i:\Gamma\to
\R$, $i=1,2$, by
\begin{equation}
\left\{ \begin{array}{ll} \overline{q}_1=-k_0-\frac{1}{6}h,\\
\overline{q}_2=k_0-\frac{1}{6}h,\\
\overline{p}_1=-\frac{1}{2}\left(k_1-k_2-k_3-3w\right),\\
\overline{p}_2=-\frac{1}{2}\left(k_1-k_2+k_3+3w\right),\\
\overline{r}_1=\frac{1}{2}\left(k_1+k_2-k_3-3w\right),\\
\overline{r}_2=-\frac{1}{2}\left(k_1+k_2+k_3-3w\right).
\end{array}\right.\label{CP5}
\end{equation}
Let us now consider the embedding
\begin{equation}\widetilde{\Gamma}=(\widehat{\mathcal{R}},\overline{q},\overline{p},\overline{r}):\Gamma\to
P.\end{equation} From (\ref{CP3}) and (\ref{CP5}) it follows that
$\widetilde{\Gamma}$ is a $1$-dimensional integral manifold of the
differential system $(\mathcal{I},\Omega)$. We set
\begin{equation}d\overline{q}_j=\overline{q}^*_j\mu,\quad
d\overline{p}_j=\overline{p}_j^*\mu,\quad
d\overline{r}_j=\overline{r}^*_j\mu,\label{CP6}\end{equation}
where $\overline{q}^*_j$, $\overline{p}^*_j$, $\overline{r}^*_j$
are real-analytic functions. From (\ref{CP3}) we infer that
\begin{equation}
\widetilde{\Gamma}_*\left(\frac{\partial}{\partial
\mu}\right)=\frac{\partial}{\partial
\omega^1}-\frac{\partial}{\partial
\omega^2}+\overline{q}^*_i\frac{\partial}{\partial
\pi^i}+\overline{p}^*_i\frac{\partial}{\partial \upsilon^i}+
\overline{r}^*_i\frac{\partial}{\partial
\zeta^i}.\label{CP7}\end{equation} Thus, $\widetilde{\Gamma}$ is a
non-characteristic $K$-regular integral curve of
$(\mathcal{I},\Omega)$. Therefore, there exist a unique
$2$-dimensional integral manifold $\widetilde{M}\subset P$ such
that $\widetilde{\Gamma}\subset \widetilde{M}$. We consider the
Legendre immersion
\begin{equation}\phi:([A],q,p,r)\in \widetilde{M}\to [A_0\wedge A_1]\in
\mathcal{K}.\label{CP8}\end{equation} Since our arguments are
local in nature, we suppose that $\phi$ is one-to-one and we
identify $\widetilde{M}$ with its image $M=\phi(\widetilde{M})$.
Then, the map
\begin{equation}\mathcal{A}:([A],q,p,r)\in M \to [A]\in \widetilde{G}\label{CP9}\end{equation} is the normal
frame field along $M$. From this we deduce that $M$ is a
Lie-minimal surface. By construction, $\Gamma$ is contained in $M$
and
\begin{equation}\mathcal{R}=\mathcal{A}|_{\Gamma},\quad
\alpha^1|_{\Gamma}=-\alpha^2|_{\Gamma}=\mu,\quad
\overline{q}_i=q_i|_{\Gamma},\quad
\overline{p}_i=p_i|_{\Gamma},\quad
\overline{r}_i=r_i|_{\Gamma},\quad
i=1,2.\label{CP10}\end{equation} In particular, the $1$-form
$\alpha^1+\alpha^2$ vanishes identically along $\Gamma$. Combining
(\ref{CP5}) and (\ref{CP10}) we deduce
$$w=\frac{1}{3}(p_1-p_2)|_{\Gamma},\quad
f=-3(q_1+q_2)|_{\Gamma},\quad \sigma_L=\sigma_0|_{\Gamma}.$$ Let
us recall that the curvature sphere mappings $\sigma_0$ and
$\sigma_1$ are represented by $[A_0]$ and by $[A_1]$ respectively.
On the other hand, $L$ is spanned by the first row vector of the
framing $\widehat{\mathcal{R}}$ so that
$\sigma_L=[\widehat{R}_0]$. This implies
$$\sigma_1|_{\Gamma}=[A_0]|_{\Gamma}=[\widehat{R}_0]=\sigma_L.$$
From this we infer that $M$ satisfies the required properties. The
uniqueness of $M$ follows from the uniqueness of the real-analytic
integral manifold $\widetilde{M}$ containing
$\widetilde{\Gamma}$.\qer}\end{proof}

\medskip

\noindent \small{\sc Dipartimento di Matematica Pura ed Applicata,
Universit\'a di L'Aquila, Via Vetoio, 67100 L'Aquila, Italy; {\tt
musso@univaq.it}

\end{document}